%% file: main.tex
\documentclass[11pt]{amsart}
\usepackage[margin=3cm]{geometry}               
\usepackage[usenames,dvipsnames,svgnames,table]{xcolor}    
\usepackage{graphicx}
\usepackage{tikz}
\usepackage{tikz-cd}
\usetikzlibrary{decorations.pathreplacing}
\usepackage{lineno, hyperref}
\usepackage{amssymb,amsfonts, amsmath, amsthm}
\usepackage{epstopdf}
\usepackage[shortlabels]{enumitem}
\usepackage{ulem}
\usepackage{cancel}
\usepackage{thmtools} 
\usepackage{mathtools} 
\usepackage{mathrsfs}
\usepackage{yhmath}
\swapnumbers
\input{commands.tex}

\title{On the density at infinity of definable functions }

\author[Si Tiep Dinh]{Si Tiep Dinh}
\address{Institute of Mathematics, VAST, 18, Hoang Quoc Viet Road, Cau Giay District 10307,
Hanoi, Vietnam}
\email{dstiep@math.ac.vn}

\author[Nhan Nguyen]{Nhan Nguyen}
\address{FPT University, Danang, Vietnam}
\email{nguyenxuanvietnhan@gmail.com}

\begin{document}

\maketitle
\begin{abstract}
In this paper, we give a simple proof that the density at infinity of fibers of a definable function is locally Lipschitz outside the set of asymptotic critical values. 
\end{abstract}

\setcounter{tocdepth}{1}


\newcounter{eqn}

\normalem

\section{Introduction}

Consider a $C^2$ definable function $f: \mathbb{R}^n \to \mathbb{R}$. Our primary focus lies in understanding the behavior of the fibers of $f$ at infinity, which can be captured by the concept of the asymptotic critical values of $f$, the set of all such values is defined as:
$$K_\infty(f) = \{ y \in \bb R: \text{ there exist a sequence } (x_l) \to \infty, f(x_l) \to y, \|x_l\|\|\nabla f (x_l) \|\to 0\}.$$

It is worth noting that $K_\infty (f)$ is a finite set containing the bifurcation set $ B_\infty(f)$ of $f$ at infinity (see \cite{D'acunto}). In particular, computing $K_\infty (f)$ is significantly easier than computing $B_\infty (f)$. Recall that the bifurcation set of $f$ at infinity is the set of all values $y$ for which $f$ is not topologically trivial at infinity.

In \cite{Grandjean_1} and \cite{Grandjean_2}, Grandjean explored the continuity of functions $t \mapsto K(t)$ and $ t \mapsto  |K|(t)$, which respectively denote the total curvature and total absolute curvature of the fiber $f^{-1}(t)$. He established that these functions admit only a finite number of discontinuities. However, the precise nature of these discontinuities had remained uncharacterized.

In a later work, Dutertre and Grandjean \cite{Nicolas_1} pointed out that the discontinuities only appear on the union of the critical values and the asymptotic critical values of $f$. 

Recently, in \cite{Tiep}, Dinh and Pham investigated the behavior of the tangent cone at infinity of the fibers of a polynomial $f$ outside the set $K_\infty(f)$. They successfully demonstrated that the density of these tangent cones is locally Lipschitz. In a different context, Dutertre and Grandjean \cite{Nicolas_2} studied the continuity of Lipschitz-Killing curvatures, a generalized notion of density, at infinity for definable functions. They established that Lipschitz-Killing curvatures at infinity of fibers of $f$ are continuous on $\mathbb{R} \setminus K_\infty(f)$.

In this paper, we present a simple  proof that the density at infinity of these fibers is locally Lipschitz. As a special case, our result strengthens the findings of Dutertre--Grandjean. The main idea of our proof is as folows: First, we establish the equivalence between the density $\theta (X, \infty)$ of a definable set at infinity and the density $\theta (\varphi (X), 0)$ of $\varphi(X)$ at the origin  where $\varphi (x) = x/\|x\|^2$. This equivalence allows us to shift the study of density at infinity to the more manageable study of density at the origin. Consider the map $H = (f(x), \varphi(x))$. It is easy to see that  the density of $f^{-1}(t)$ at infinity is equal to the density of the fiber of $H(X)$ at $t$, with respect to the projection $\pi: \bb R \times \bb R^n \to \bb R$. We show that for each connected component $I$ of $\bb R\setminus K_{\infty}(f)$, the collection $\{H(f^{-1}(I)), I\times{0}^n\}$ satisfies the $(w)$-regularity condition. Then, we apply \cite[Proposition 4.4]{nhan} to arrive at the desired result.

Throughout the paper, we assume the reader's familiarity with the notion of o-minimal structures on $\mathbb{R}$. For more comprehensive details,  we refer the reader to \cite{Coste}, \cite{Dries}, \cite{Loi1}.

We denote $\mathbb{B}^{n}_r$ and $\mathbb{S}^{n-1}_r$ respectively the $n$-dimensional closed ball and the $(n-1)$-dimensional sphere in $\mathbb{R}^n$ of radius $r$ centered at $0$. For a subset $X \subset \mathbb{R}^n$, $\overline{X}$ denotes the closure of $X$ in $\mathbb{R}^n$. Given non-negative functions $f, g: X \to \mathbb{R}$, we use write $f\lesssim g$ if there exists a constant $C>0$ such that $f(x) \leq C g(x)$ for all $x \in X$. Such a constant $C$ is referred to as {\it a constant for the relation $\lesssim$}.

\subsection{Acknowledgements}

The authors would like to express their gratitude to the Vietnam Institute for Advanced Study in Mathematics (VIASM) for their warm hospitality and generous support during the writing of this paper.

\section{Existence of density at infinity}
Let $X$ be a definable set in $\Bbb R^n$ and let $x \in \overline{X}$. Suppose that $\dim (X,x) = k$, i.e., $\dim (X \cap \bb B_r^n(x)) = k$ for $r>0$ sufficiently small. The \textbf{density} of $X$ at $x$ is defined as 
$$\displaystyle\theta (X,x):=\lim_{r\to 0}\frac{\vol_k\big(X\cap\Bbb B^n_r(x)\big)}{\vol_k(\Bbb B^k_r)}. $$
It has been known that $\theta (X, x)$ always exists (see for example \cite{Kur}, \cite{Comte}). 

In the following, we show that the $k-$dimensional density of $X$ at $x$ can also be defined by using the volume of the intersection of $X$ with spheres centered at $x$ instead of balls.
\begin{prop}\label{SphericalPro} We have 
\begin{equation}\label{Spherical}
\displaystyle\theta (X,x)=\lim_{r\to 0}\frac{\vol_{k-1}\big(X\cap\Bbb S^{n-1}_r(x)\big)}{\vol_{k-1}(\Bbb S^{k-1}_r)}.
\end{equation}
\end{prop}
\begin{proof} 

Without loss of generality, we assume that $X$ is a $C^1$ manifold. If it were not the case, we could decompose $X$ into a finite union of $C^1$ definable submanifolds. Then, we would consider the union of those submanifolds that have the highest dimension and contain $x$ in their adherence. We endow $X$ with the Riemannian structure induced by the Euclidean structure on $\mathbb{R}^n$. Additionally, we identify $x$ with the origin $0$.

Let $\rho:\mathbb{R}^n\to\mathbb{R}$ be defined as $\rho(y)=|y|$. Since the set of critical values of the restriction $\rho|_X$ is finite, we can find an $\epsilon>0$ such that $\rho|_X$ does not have critical values in the interval $(0,\epsilon)$. Consequently, $X\cap\mathbb{S}^{n-1}_r$ becomes an $C^1$ manifold.

We further endow $X\cap\mathbb{S}^{n-1}_r$ with the Riemannian structure induced by the Riemannian structure on $X$, which coincides with the Riemannian structure induced by the Euclidean structure on $\mathbb{R}^n$. The corresponding volume density is denoted by $dV{\rho|_X^{-1}(r)}(y)$.  Applying the co-area formula (see e.g.,  \cite[formula (9.1.7)]{Nicolaescu_2}), we obtain: 
$$\vol_k\big(X\cap\Bbb B^{n}_r \big)=\int_{t=0}^r\Big(\int_{\rho|_X=t}\frac{1}{\|\nabla(\rho|_X)(y)\|}\Big|dV_{\rho|_X^{-1}(t)}(y)\Big|\Big)dt.$$
Moreover, it is clear that 
$$\vol_k\big(\Bbb B^{k}_r\big)=\int_{t=0}^r\Big(\int_{\Bbb S^{k-1}_t}\Big|dV_{\Bbb S^{k-1}_t}(y)\Big|\Big)dt=\int_{t=0}^r\vol_{k-1}(\Bbb S^{k-1}_t)dt.$$
Hence 
$$\displaystyle\theta (X,0)=\lim_{r\to 0}\frac{\displaystyle\int_{t=0}^r\Big(\int_{\rho|_X=t}\frac{1}{\|\nabla(\rho|_X)(y)\|}\Big|dV_{\rho|_X^{-1}(t)}(y)\Big|\Big)dt}{\displaystyle\int_{t=0}^r\vol_{k-1}(\Bbb S^{k-1}_t)dt}.$$
By l'H\^opital's rule,  if the limit 
\begin{equation}\label{equ3}
    \lim_{r\to 0}\frac{\displaystyle\int_{\rho|_X=r}\frac{1}{\|\nabla(\rho|_X)(y)\|}\Big|dV_{\rho|_X^{-1}(r)}(y)\Big|dr}{\vol_{k-1}(\Bbb S^{k-1}_r)}
\end{equation}
exists then it is equal to $\theta (X,0)$.

Since $\displaystyle\int_{\rho|_X=r}\Big|dV_{\rho|_X^{-1}(r)}(y)\Big|=\displaystyle\vol_{k-1}\big(X\cap\Bbb S^{n-1}_r(x)\big)$, to prove the proposition, it suffices to show that 
$$\|\nabla(\rho|_X)(y)\|\to 1 \text{ as } y \to 0.$$

Note that $\nabla(\rho|_X)(y)=\proj_{T_yX}\nabla\rho(y)=\displaystyle\proj_{T_yX}\Big(\frac{y}{\|y\|}\Big),$ 
where $\proj_{T_yX}$ denotes  the projection onto $T_yX$. Since $\|\nabla\rho(y)\|=1$, it follows that $\nabla(\rho|_X)(y)\leqslant 1$. Assume on the contradictory that  that there is a sequence $y^k\to 0$ such that $\|\nabla(\rho|_X)(y^k)\|\to c<1.$ By Curve Selection Lemma, there is a $C^1$ definable curve  $\varphi:(0,\epsilon)\to X$ such that $\lim_{r\to 0}\varphi(r)=0$ and $\lim_{r\to 0}\|\nabla(\rho|_X)\big(\varphi(r)\big)\|\to c$. Reparametrizing $\varphi$ we may assume  $\|\varphi(r)\|=r$. Since $\varphi$ is definable, $\lim_{r\to 0}\frac{\varphi(r)}{\|\varphi(r)\|}$ exists. 

By l'H\^opital's rule: 
$$\lim_{r\to 0} \frac{\varphi(r)}{\|\varphi(r)\|} = \lim_{r\to 0} \frac{\varphi'(r)\|\varphi(r)\|}{\langle \varphi(r), \varphi'(r)\rangle} = \lim_{r\to 0} \frac{\varphi'(r)}{\|\varphi'(r)\| \cos(\widehat{\varphi(r),\varphi'(r)})} = \lim_{r\to 0} \frac{\varphi'(r)}{\|\varphi'(r)\|}$$
since $(\widehat{\varphi(r),\varphi'(r)})\to 0$ as $r \to 0$.

Consequently
$$\lim_{r\to 0}\|\nabla(\rho|_X)\big(\varphi(r)\big)\|=\lim_{r\to 0}\Big\|\proj_{T_{\varphi(r)}X}\Big(\frac{\varphi(r)}{\|\varphi(r)\|}\Big)\Big\|=\lim_{r\to 0}\Big\|\proj_{T_{\varphi(r)}X}\Big(\frac{\varphi'(r)}{\|\varphi'(r)\|}\Big)\Big\|.$$
Since $\varphi(r)\in X,$ it follows that $\varphi'(r)\in T_{\varphi(r)}X.$ Thus $\displaystyle\proj_{T_{\varphi(r)}X}\Big(\frac{\varphi'(r)}{\|\varphi'(r)\|}\Big)=\frac{\varphi'(r)}{\|\varphi'(r)\|}$ and we deduce that
$\lim_{r\to 0}\|\nabla(\rho|_X)\big(\varphi(r)\big)\|=1,$ which is a contradiction.
\end{proof}

 Now we suppose that $X$ is unbounded, and $\dim (X, \infty) = k$, i.e., $\dim(X\setminus\Bbb B^n_r)=k$ for $k$ sufficiently large. The following result enables us to define the density at infinity for $X$ and establishes a connection between the density at infinity and the density of a suitable set at the origin.
 
\begin{thm}\label{InfinityLocal} Let $X$ be a definable set of dimension $k$ at infinity. Then the limit $$\displaystyle\theta (X,\infty):=\lim_{r\to +\infty}\frac{\vol_{k-1}(X\cap\Bbb S^{n-1}_r)}{\vol_{k-1}(\Bbb S^{k-1}_r)}$$ exists, moreover, 
$$\theta (X,\infty)=\theta (\varphi(X\setminus\{0\}),0),$$
where $\varphi:\R^n\setminus\{0\}\to\R^n\setminus\{0\}$ is the map given by $\varphi(x):=\displaystyle\frac{x}{\|x\|^2}.$
\end{thm}
\begin{proof} Without any loss of generality, let's assume that $0\not\in X$. To prove the theorem, it suffices to show that $\theta(X,\infty)=\theta(\varphi(X),0)$. Now, we define $X_r:=X\cap\mathbb{S}^{n-1}_r$. The following lemma establishes a relationship between the $(k-1)$-dimensional volume of $X_r$ and $\varphi(X_r).$
\begin{lem} \label{InverseVolume} We have $\displaystyle\vol_{k-1}\big(\varphi(X_r)\big)=\frac{1}{r^{2(k-1)}}\vol_{k-1}(X_r).$
\end{lem}
\begin{proof}
 By definition of $(k-1)-$dimensional Hausdorff measure,
$\vol_{k-1}(X_r)=\displaystyle\lim_{\delta\to 0}H^{k-1}_\delta(X_r)$ where
$$H^{k-1}_\delta(X_r)=\inf\Big\{\sum_{i=1}^\infty(\diam U_i)^{k-1}:\ U_i\subset\Bbb S^{n-1}_r,\ \bigcup_{i=1}^\infty U_i\supset X_r, \diam U_i<\delta\Big\},$$
$\diam U_i$ is the geodesic diameter of $U_i$ with respect to the induced Riemannian structure on $\Bbb S^{n-1}_r$. 
Let us prove that 
\begin{equation}\label{H}
\displaystyle H^{k-1}_{\frac{\delta}{r^2}}\big(\varphi(X_r)\big)=\frac{1}{r^{2(k-1)}}H^{k-1}_\delta(X_r).
\end{equation}
Let $V_i:=\varphi(U_i)$, then it is easily seen that $V_i\subset\Bbb S^{n-1}_{\frac{1}{r}}$ and $\displaystyle\bigcup_{i=1}^\infty V_i\supset \varphi(X_r)$. Moreover, it is not hard to check that $\displaystyle\diam V_i=\frac{1}{r^2}\diam U_i$, so $\displaystyle\diam V_i<\frac{\delta}{r^2}$. Hence, we  
\begin{equation}\label{H1}
\displaystyle H^{k-1}_{\frac{\delta}{r^2}}\big(\varphi(X_r)\big)\leqslant\frac{1}{r^{2(k-1)}}H^{k-1}_\delta(X_r).
\end{equation}
Applying Inequality (\ref{H1}) with $X_r$ replaced by $\varphi(X_r)\subset \Bbb S^{n-1}_{\frac{1}{r}}$, we obtain
$$\displaystyle H^{k-1}_{r^2{\delta}}(X_r)=H^{k-1}_{r^2{\delta}}\Big(\varphi\big(\varphi(X_r)\big)\Big)\leqslant{r^{2(k-1)}}H^{k-1}_\delta\big(\varphi(X_r)\big).$$
Replacing $\delta$ by $\displaystyle\frac{\delta}{r^2}$, we get
\begin{equation}\label{H2}
\displaystyle H^{k-1}_{\delta}(X_r)\leqslant{r^{2(k-1)}}H^{k-1}_{\frac{\delta}{r^2}}\big(\varphi(X_r)\big).
\end{equation}
Now, (\ref{H}) follows (\ref{H1}) and (\ref{H2}). Then the lemma is proved by letting $\delta\to 0.$
\end{proof}
By Lemma \ref{InverseVolume}, we have
$$\vol_{k-1}(X_r)=r^{2(k-1)}\vol_{k-1}\big(\varphi(X_r)\big)=r^{2(k-1)}\vol_{k-1}\big(\varphi(X)\cap\Bbb S^{k-1}_{\frac{1}{r}}\big).$$
This, together with Proposition \ref{SphericalPro}, yields
$$\begin{array}{llll}
\displaystyle\theta (X,\infty)&=&\displaystyle\lim_{r\to +\infty}\frac{\vol_{k-1}(X_r)}{\vol_{k-1}(\Bbb S^{k-1}_r)}\\
&=&\displaystyle\lim_{r\to +\infty}\frac{r^{2(k-1)}\vol_{k-1}\big(\varphi(X)\cap\Bbb S^{k-1}_{\frac{1}{r}}\big)}{r^{2(k-1)}\vol_{k-1}(\Bbb S^{k-1}_{\frac{1}{r}})}\\
&=&\displaystyle\lim_{r\to +\infty}\frac{\vol_{k-1}\big(\varphi(X)\cap\Bbb S^{k-1}_{\frac{1}{r}}\big)}{\vol_{k-1}(\Bbb S^{k-1}_{\frac{1}{r}})}\\
&=&\theta(\varphi(X),0).
\end{array}$$
The theorem is proved.
\end{proof}

\begin{rem}\rm
The limit $\theta (X, \infty)$ is referred as \textbf{the density of $X$ at infinity}. In the case that $\dim X = \dim (X, \infty)$, it is easily shown  that the limit $\displaystyle\lim_{r\to +\infty}\frac{\vol_k(X\cap\Bbb B^n_r)}{\vol_k(\Bbb B^k_r)}$ is equal to $\theta (X,\infty)$ by utilizing the co-area formula. The detailed proof  is  left to the reader. 
Consequently, we have
\begin{equation}\label{DensityInfinityDef}\displaystyle\theta (X,\infty)=\displaystyle\lim_{r\to +\infty}\frac{\vol_k(X\cap\Bbb B^n_r)}{\vol_k(\Bbb B^k_r)}=\lim_{r\to +\infty}\frac{\vol_{k-1}(X\cap\Bbb S^{n-1}_r)}{\vol_{k-1}(\Bbb S^{k-1}_r)}=\theta (\varphi(X\setminus\{0\}),0).
\end{equation}
\end{rem}

\section{Continuity of density at infinity}
 
Let $f: \bb R^n \to \bb R$ be a $C^2$ definable function ($n \geq 2$). Our main result is as follows:

\begin{thm} \label{thm_main}
The function $t \mapsto   \theta(f^{-1}(t), \infty)$ is locally Lipschitz on $\bb R \setminus K_\infty (f)$.
\end{thm} 

We need the following lemmas to prove the above theorem.

\begin{lem}\label{lem3.1} Let $I = (a, b)$ be an interval in $\mathbb{R}$, and let $g: I \to \mathbb{R}$ be a positive definable function satisfying the following property: for every $t \in I$ and for every sequence $(t_n)_{n\in\mathbb{N}} \subset I$ tending to $t$, if the limit $\lim_{n\to \infty} g(t_n)$ exists, then it is positive.

Then, there exists a positive definable continuous function $h: I \to \mathbb{R}$ such that $h(t) \leq g(t)$ for all $t \in I$.
\end{lem}

\begin{proof}
    Since $g$ is definable, there exist finite points $a = a_1 < a_2 < \ldots < a_k = b$ such that $g$ is continuous on each open interval $(a_i, a_{i+1})$, where $i = 0, \ldots, k-1$. For each $0 < i < k$, we set:
 $$m_i = \lim_{t \to a_i^-} g(t)  \text{ and } n_i = \lim_{t \to a_i^+} g(t).$$
From the given assumption, it follows that $m_i$ and $n_i$ are positive. Now, for each $0 \leq i < k$, we let $g_i: I \to \mathbb{R}$ be a continuous extension of $g|_{(a_i, a_{i+1})}$ to $I$ by letting it be $n_i$ on $(a, a_i]$ if $i > 0$, and $m_{i+1}$ on $[a_{i+1}, b)$ if $i < k$.

It is clear that $g_i$ are definable positive continuous functions. We then define $h: I \to \mathbb{R}$ as $h = \min_{0\leq i \leq k}{g_i}$. Hence, $h$ is the desired function.
  
\end{proof}

Employing similar arguments as in the proof of Lemma \ref{lem3.1}, we obtain the following lemma.
\begin{lem}\label{lem3.2} Let $I = (a, b)$ be an interval in $\bb R$ and $g: I \to \bb R$ be a non-negative definable function such that for every  $t\in I$ and for every sequence $(t_n)_{n\in\bb N}\subset I$ tending to $t$ we have $\lim_{n\to \infty} g(t_n) < \infty$. Then, there is a positive continuous definable function $h: I \to \bb R$ such that $h(t)\geq g(t)$ for all $t \in I$.
\end{lem}

\begin{proof}[Proof of Theorem \ref{thm_main}] Since $K_\infty (f)$ is finite,  $\bb R \setminus K_\infty (f)$ is a disjoint union of open intervals. Let $I = (a, b)$ be such an interval.  It suffices to prove that the function  $t\mapsto \theta(f^{-1}(t), \infty))$ is locally Lipschitz on $I$.  Note that the restriction $f|_{f^{-1} (I)}: f^{-1} (I) \to I$ is topologically trivial at infinity. Without loss of generality, we may assume that $f^{-1}(t)$ is unbounded for all $t\in I$, since otherwise the germ at infinity of $f^{-1}(t)$ is empty, hence the density is zero, and the result is trivial. 

\begin{claim}\label{claim_1}
   There are  continuous  definable functions  $\varepsilon_1, \varepsilon_2: I \to \bb R_+$ such that 
$$ \forall t\in I,\forall  x \in f^{-1}(t), \|x\| > \varepsilon_1(t) \Rightarrow \|x\|\|\nabla f(x)\| > \varepsilon_2(t).$$
\end{claim}

\begin{proof}
Let us give a proof of Claim \ref{claim_1}. For $t \in I$, we define 
$$
\sigma_1(t) :=
\begin{cases}
1, \text{ if } \nabla f(x) \neq 0 \text{ for all } x \in f^{-1}(t)\\
 1 + \sup\{\|x\|: x \in f^{-1}(t) \text{ and } \nabla f (x) = 0\}  \text{ otherwise}.
\end{cases}
$$
Since $t \not\in K_{\infty}(f)$, it is easy to check that for every $t\in I$,  $1\leq \sigma_1(t) < \infty$ and moreover, if $ (t_n) \subset I$ is a sequence tending to $t$ then $\lim_{n\to \infty} \sigma_1(t_n) < \infty$. Applying Lemma  \ref{lem3.2} we obtain a continuous definable function  $\varepsilon_1: I \to \bb R_+$ such that $\varepsilon_1 > \sigma_1$. 
For $t \in I$, set
$$\sigma_2(t)  := \inf \{\|x \| \|\nabla f(x)\|: f(x) = t \text { and } \| x \|\geq \varepsilon_1(t)\}.$$
Again, since $t \not\in K_{\infty}(f)$ we see that $\sigma_2(t)$ is a positive definable function satisfying the condition in Lemma \ref{lem3.1}. Choose  $\varepsilon_2: I \to \bb R_+$ to be a continuous definable such that $\varepsilon_2 < \sigma_2$. Then, $\varepsilon_1$ and $\varepsilon_2$ are functions satisfying Claim 1. 
\end{proof}

Now, let us define $X := \{ x \in f^{-1}(I): \|x\| > \varepsilon_1(f(x))\}$, which is an open subset of $\mathbb{R}^n$. Consider the map $H: \mathbb{R}^n \to \mathbb{R} \times \mathbb{R}^n$ defined by $ H(x): = (f(x), \frac{x}{\|x\|^2})$. Observe that $H$ is an injection which is of class $C^2$ everywhere except at the origin. Additionally, $\text{rank} d_xH = n$ for $x \neq 0$. This implies that $H|_{\mathbb{R}^n\setminus {0}}$ is a $C^2$ embedding. Consequently, $Z:= H(X)$ forms a $C^2$ submanifold of $\mathbb{R} \times \mathbb{R}^n$.

Let $Y := I\times {0}^n$. It is obvious that $Y \subset \overline{Z}$. Thus, $(Z, Y)$ forms a stratification.

\begin{claim}\label{claim_2}
$(Z, Y)$ is $(w)$-regular.
\end{claim}

\begin{proof}
To prove the claim, it suffices to show that the vector field $\frac{\partial }{\partial t}$ defined on $Y$ can be extended to a rugose stratified vector field $v(t, u)$ on $Z \cup Y$ (see \cite[Proposition 2]{BT}). 
In fact, on $Z$ the extension is defined to be the pushforward of the vector field $\frac{\nabla f(x)}{\|\nabla f(x)\|^2}$ under the map $dH$. More precisely, 

$$
v(t,u): =
\begin{cases}
\partial /\partial t, \text{ for } (t,u) \in Y\\
dH_x \displaystyle \left(\frac{\nabla f(x)}{\|\nabla f(x)\|^2} \right) =  \frac{\partial }{\partial t} + d \varphi_x  \left( \frac{\nabla f(x)}{\|\nabla f(x)\|^2} \right) , \text{ for } (t,u) \in Z\\
\end{cases}
$$
 where $x := H^{-1}(t, u)= \frac{u}{\|u\|^2}$ and $\varphi(x): =\frac{x}{\|x\|^2}$.

It remains to check the rugosity of $v$, which means for each $a = (t, 0) \in Y$, there exists a constant $C>0$ and a neighborhood $U_a$ of $a$ in $\mathbb{R} \times \mathbb{R}^n$ such that
$$ \|v(z) - v(y)\| \leq C \|z - y\|$$
for all $z \in Z \cap U_a$ and $y \in Y \cap U_a$.

Given $a_0 = (t_0,0) \in Y$, choose a small neighborhood $U$ of $a_0$ in $\bb R^{n+1}$. Since $\varepsilon_2 (t)$ obtained from Claim \ref{claim_1} is continuous on $I$, there is a constant $c >0$ such that $\varepsilon_2(t) > c$ for every $a = (t,0) \in U \cap Y$.  Let $z = (t_z, u_z) \in Z \cap U$ and $y = (t_y, 0) \in Y \cap U$. We have 
$$\|v(z) - v(y)\| \leq \|v(z) - v(t_z, 0)\| +\| v(t_z, 0) - v(t_y, 0)\| = \|v(z) - v(t_z, 0)\| = \left\|d_x\varphi  \left(\frac{\nabla f(x)}{\|\nabla f(x)\|^2}\right)\right\|.$$

Note that
$$
   d_x \varphi  =
  \left[ {\begin{array}{cccc}
   \displaystyle \frac{\|x\|^2 - 2x_1^2}{\|x\|^4}  &  \displaystyle \frac{-2x_1 x_2}{\|x\|^4} &  \dots  &  \displaystyle  \frac{-2x_1 x_n}{\|x\|^4}\\
   
   \vdots & \vdots&  & \vdots \\
   
   \displaystyle \frac{- 2x_n x_1}{\|x\|^4}   &  \displaystyle \frac{-2x_n x_2}{\|x\|^4} &  \dots  &   \displaystyle \frac{\|x\|^2 - 2x_n^2}{\|x\|^4}
     \end{array} } \right]
$$
We see that the absolute value of each entry of $d_x\varphi $ is $\lesssim  \frac{1}{\|x\|^2}$ so $\|d_x\varphi\| \lesssim \frac{1}{\|x\|^2}$. This yields that  
$$\left\|d_x\varphi\left(\frac{\nabla f(x)}{\|\nabla f(x)\|^2}\right)\right\| \lesssim  \|d_x\varphi \| \bigg \|\frac{\nabla f(x)}{\|\nabla f(x)\|^2}\bigg\| \lesssim \frac{1}{\|x\|^2\|\nabla f(x)\|} \lesssim  \frac{1}{\varepsilon_2(f(x))\|x\|} <  \frac{1}{c\|x\|}  .$$
On the other hand, 
$$ \|z - y\|  \geq \| z - (t_z, 0)\| = \| u_z\| = \frac{1}{\|x\|}.$$
This implies that 
$$ \|v(z) -v( y)\| \lesssim \|z - y\|.$$
Note that the constant for the relation $\lesssim$ depends only on $U_a$. Thus, Claim \ref{claim_2} is proved.  
\end{proof}

Now applying \cite[Proposition 4.4]{nhan}, we get $\theta (Z_t, 0)$ is locally Lipschitz in $t$ where $Z_t := \{u\in \bb R^n: (t, u) \in Z\}$. By Theorem \ref{InfinityLocal},  $\theta(Z_t, 0) = \theta (f^{-1}(t), \infty)$, then the theorem follows. 
\end{proof}

\bibliographystyle{siam}
\bibliography{Biblio.bib}

\end{document}

%% file: Commands.tex
\newtheorem*{theorem*}{Theorem}
\newtheorem*{question*}{Question}
\newtheorem*{corollary*}{Corollary}

\numberwithin{equation}{subsection}

\newtheorem{thm}{Theorem}[section]
\newtheorem{lem}[thm]{Lemma}
\newtheorem{prop}[thm]{Proposition}

\newtheorem{rem}[thm]{Remark}

\newtheorem{claim}[thm]{Claim}

\numberwithin{equation}{section}
\newcommand{\bb}{\mathbb}

\newcommand{\vol}{{\rm vol}}

\newcommand{\proj}{{\rm proj}}
\newcommand{\diam}{{\rm diam}}

\newcommand{\R}{\Bbb R}